# Hamiltonian formalism for optimal control of nonlinear loaded integro-PDEs.

S. A. Belbas

Independent consultant; formerly, with the Dept. of Mathematics, University of Alabama, Tuscaloosa, AL 35487-0350, USA. e-mail SBELBAS@GMAIL.COM

**Abstract**

We formulate nonlinear nonlocal integro-PDE with memory, biloaded (boundary integrals load the ambient space, and the ambient space loads the boundary), and the associated optimal control problems. We derive part of the necessary conditions for optimality in the form of Hamilton-Euler-Lagrange loaded integro-PDEs. In the process, we introduce an agglomeration of new differential operators.
Our results have relevance to optimal amelioration of flooded areas, remediation of sites of contaminated groundwater, and active control methods for optimally extinguishing forest fires.

## 1 Introduction

We deal with the derivation of the Hamiltonian equations for certain systems of nonlinear integro-PDEs of three main types, classified according to the differential order (diff-order): first-order, second-order, and Barenblatt [4, 5] third order (first order with respect to time, and second order with respect to the spatial coordinates). For integro-differential equations, there is also the possibility of order of integral multiplicity (int-order) if the integro-differential equations include multiple integrations over the ambient space, e.g. [23]; we do not consider integral multiplicity higher than 1 in this work. The three types of integro-differential equations correspond to 3 physical situations: control methods for the amelioration of flooded areas (first diff-order integro-PDE, shallow water waves [11, 12, 13, 14), remediation of contaminated groundwater sites [15, 16], and control-theoretic methods for extinguishing forest fires [17, 18]. The paper belongs to the genre "cybernetic physics", for which genre the archetypical reference, to our perception, is [1], and closely following are [2, 3]. Relevant essential works also include [7, 8] (nonlinear and nonlocal, resp., thermodynamics). We use the framework of third-kind integral equations with partial derivatives, and loaded equations [9,10], which are also reversely loaded (the equations over the ambient space include integrals over lower-dimensional manifolds, and the equations for the boundary terms include also integrals over the ambient space). Existence theorems are not included in this work, and they are available only for particular cases, e.g. [6] for Barenblatt-type PDEs in one space dimension. Hamiltonian equations are instrumental in synthesizing optimal policies for control of integral systems [22]. The rigorous theory of necessary optimality conditions for integral systems originates with [21].

## 2 The Barenblatt third-diff-order loaded systems

We shall formulate the controlled integro-PDE of Barenblatt third-differential order, from which the Hamiltonian equations of all 3 types of integro-PDE loaded systems of interest for optimal control problems with relevance to environmental sciences, will follow. We use the formulation of third-kind integral equations. The spatiotemporal ambient domain is $Q = \Omega \times (0,T) \subseteq \mathbb{R}^d \times \mathbb{R}$, the spatial domain $\Omega$ is an open bounded set in $d$ – dimensional Euclidean space with sufficiently smooth boundary $\Gamma = \partial\Omega$. The hypersurface area element on $\Gamma$ will be denoted by $d\sigma$. The state (in the sense of control theory) will be an $n$ – dimensional vector-valued function $\varphi$. The spatial gradient of the state, i.e. matrix with $ij$ – th element $\dfrac{\partial \varphi_i}{\partial x_j}$, will be $p = [p_i^j]$. A dash indicates time-derivative. The 3-index array of the second derivatives is $q = [q_i^{jk}] = \left[\dfrac{\partial^2 \varphi_i}{\partial x_j \partial x_k}\right]$. Greek letters, in the arguments of a function, indicate traces onto the boundary. $u$ will be the control in $Q$, and $w$ the control on $(0, T) \times \Gamma$. Subscripts 0 and $T$, in the controls and the state, will indicate controls and state at the initial and the final time, respectively.

We shall use the following notation:

$$S(t,x) = (\varphi, p, q, \varphi', p', q', u)(t,x); \; S_\partial(t,\xi) = (\varphi, \varphi', p, p', w)(t,\xi); \; S_0(x) = (\varphi_0, p_0, q_0, u_0)(x); \; S_T(x) = (\varphi_T, p_T, q_T, u_T)(x);$$
$$S_{0\partial}(\xi) = (\varphi_0, p_0, w_0)(\xi); \; S_{T\partial}(\xi) = (\varphi_T, p_T, w_T)(\xi)$$

$$\varphi(t,x) = f_0(t,x,S_\Omega(t,x)) + \int_0^t f_1(t,x,s,S_\Omega(s,x))\,ds + \int_\Omega f_2(t,x,y,S_\Omega(t,y))\,dy +$$
$$+ \int_0^t \int_\Omega f_3(t,x,s,y,S_\Omega(s,y))\,dy\,ds + \int_\Gamma f_4(t,x,\eta,S_\Gamma(t,\eta))\,d\sigma(\eta) + \int_0^t \int_\Gamma f_5(t,x,s,\eta,S_\Gamma(s,\eta))\,d\sigma(\eta)\,ds;$$

$$\varphi(t,\xi) = g_0(t,\xi,S_\Gamma(t,\xi)) + \int_0^t g_1(t,\xi,s,S_\Gamma(s,\xi))\,ds + \int_\Omega g_2(t,\xi,y,S_\Omega(t,y))\,dy +$$
$$+ \int_0^t \int_\Omega g_3(t,\xi,s,y,S_\Omega(s,y))\,dy\,ds + \int_\Gamma g_4(t,\xi,\eta,S_\Gamma(t,\eta))\,d\sigma(\eta) + \int_0^t \int_\Gamma g_5(t,\xi,s,\eta,S_\Gamma(s,\eta))\,d\sigma(\eta)\,ds$$

$$\varphi_0(x) = f_{00}(x,S_0(x)) + \int_\Omega f_{02}(x,y,S_0(y))\,dy + \int_\Gamma f_{04}(x,\eta,S_{0\partial}(\eta))\,d\sigma(\eta);$$

$$\varphi_T(x) = f_{T0}(x,S_T(x)) + \int_0^T f_{T1}(x,s,S(s,x))\,ds + \int_\Omega f_{T2}(x,y,S_T(y))\,dy + \int_0^T \int_\Omega f_{T3}(x,s,y,S(s,y))\,dy\,ds +$$
$$+ \int_\Gamma f_{T4}(x,\eta,S_{T\partial}(\eta))\,d\sigma(\eta) + \int_0^T \int_\Gamma f_{T5}(x,s,\eta,S_\partial(s,\eta))\,d\sigma(\eta);$$

$$\varphi_0(\xi) = g_{00}(\xi,S_{0\partial}(\xi)) + \int_\Omega g_{02}(\xi,y,S_0(y))\,dy + \int_\Gamma g_{04}(x,\eta,S_{0\partial}(\eta))\,d\sigma(\eta);$$

$$\varphi_T(\xi) = g_{T0}(\xi,S_{T\partial}(\xi)) + \int_0^T g_{T1}(\xi,s,S_\partial(s,\xi))\,ds + \int_\Omega g_{T2}(\xi,y,S_T(y))\,dy + \int_0^T \int_\Omega g_{T3}(\xi,s,y,S(s,y))\,dy\,ds +$$
$$+ \int_\Gamma g_{T4}(\xi,\eta,S_{T\partial}(\eta))\,d\sigma(\eta) + \int_0^T \int_\Gamma g_{T5}(\xi,s,\eta,S_\partial(s,\eta))\,d\sigma(\eta)\,ds$$

In order to keep the problem amenable to Hamiltonian treatment, the components $\varphi_0, \varphi_T, \varphi_{0\partial}, \varphi_{T\partial}$ must be introduced as separate components of the state, and the associated controls as separate and independent controls. The price for "Hamiltonianity" is jumps at $t = 0$ and $t = T$, in general $\varphi(0^+, x) \neq \varphi_0(x)$, and so on for the remaining subscripted components of the state.

An optimal control problem concerns the minimization of a cost functional

$$J = \int_\Omega F_0(x,S_0(x),S_T(x))\,dx + \int_\Gamma G_0(\xi,S_{0\partial}(\xi),S_{T\partial}(\xi))\,d\sigma(\xi) +$$
$$+ \int_0^T \int_\Omega F_1(t,x,S(t,x))\,dx\,dt + \int_0^T \int_\Gamma G_1(t,\xi,S_\partial(t,\xi))\,d\sigma(\xi)\,dt$$

We shall need (among other things) a second-order Gauss-Ostrogradsky divergence theorem. For sufficiently smooth vector-valued functions $\varphi, \psi$ of compatible dimensions, we denote by $D_\xi \varphi$ the trace of the spatial gradient of $\varphi$ onto $\Gamma$, and by $D_\xi^*$ the skew-adjoint of $D_\xi \varphi$, i.e. $\int_\Gamma \left(\psi D_\xi \varphi + (D_\xi^* \psi)\varphi\right) d\sigma = 0$. The existence of the skew-adjoint follows from Stokes' theorem on $\Gamma$, and from $\Gamma$ having empty boundary (here, "boundary" is understood in the context of simplicial chains).

## 3. Hamiltonian and the Hamilton-Euler-Lagrange loaded integro-PDEs

The Hamiltonian <u>functional</u> contains 6 co-states (mathematically analogous to generalized momenta of classical Hamiltonian Dynamics) $\psi, \psi_0, \psi_T, \omega, \omega_0, \omega_T$ co-vector valued functions, corresponding to the state components $\varphi, \varphi_0, \varphi_T, \varphi_\partial, \varphi_{0\partial}, \varphi_{T\partial}$, and it is defined as

$$H = F_0(x,\ldots) + G_0(\xi,\ldots) + F_1(t,x,\ldots) + G_1(t,\xi,\ldots) + \psi(t,x)f_0(t,x,\ldots) + \int_t^T \psi(s,x) f_1(s,x,t,\ldots)ds +$$

$$+ \int_\Omega \psi(t,y)[f_2(t,y,x,\ldots) + f_4(t,y,\xi,\ldots)]dy + \int_t^T \int_\Omega \psi(s,y)[f_3(s,y,t,x,\ldots) + f_5(s,y,t,\xi,\ldots)]dyds +$$

$$+ \omega(t,\xi)g_0(t,\xi,\ldots) + \int_t^T \omega(s,\xi)g_1(s,\xi,t,\ldots)ds + \int_\Gamma \omega(t,\eta)[g_2(t,\eta,x,\ldots) + g_4(t,\eta,\xi,\ldots)]d\sigma(\eta) +$$

$$+ \int_t^T \int_\Gamma \omega(s,\eta)[g_3(s,\eta,t,x,\ldots) + g_5(s,\eta,t,\xi,\ldots)]d\sigma(\eta)ds + \psi_0(x)f_{00}(x,\ldots) + \int_\Omega \psi_0(y)[f_{02}(y,x,\ldots) +$$

$$+ f_{04}(y,\xi,\ldots)]dy + \omega_0(\xi)g_{00}(\xi,\ldots) + \int_\Gamma \omega_0(\eta)[g_{02}(\eta,x,\ldots) + g_{04}(\eta,\xi,\ldots)]d\sigma(\eta) +$$

$$+ \psi_T(x)[f_{T0}(x,\ldots) + f_{T1}(x,t,\ldots)] + \int_\Omega \psi_T(y)[f_{T2}(y,x,\ldots) + f_{T3}(y,t,x,\ldots) + f_{T4}(y,\xi,\ldots) + f_{T5}(y,t,\xi,\ldots)]dy +$$

$$+ \omega_T(\xi)[g_{T0}(\xi,\ldots) + g_{T1}(\xi,t,\ldots)] + \int_\Gamma \omega_T(\eta)[g_{T2}(\eta,x,\ldots) + g_{T3}(\eta,t,x,\ldots) + g_{T4}(\eta,\xi,\ldots) + g_{T5}(\eta,t,\xi,\ldots)]d\sigma(\eta)$$

The differential operators needed for the Hamiltonian equations will be easier to describe if we introduce first some notation about operations with single-index and double index arrays.
If

$$a = (a_i : 1 \leq i \leq m), \ b = (b_j : 1 \leq j \leq n), \ c = (c_i : 1 \leq i \leq m),$$
$$A = (A_{ij} : 1 \leq i \leq m, 1 \leq j \leq n), \ B = (B_{ij} : 1 \leq i \leq m, 1 \leq j \leq n),$$

then we define

$$\langle a, c \rangle = \sum_i a_i c_i, \ \langle A, B \rangle = \sum_{i,j} A_{ij} B_{ij}; \ \langle a | A | = ((aA)_j : 1 \leq j \leq n), \ (aA)_j = \sum_i a_i A_{ij};$$

$$|A|b\rangle = ((Ab)_i : 1 \leq i \leq m), \ (Ab)_i = \sum_j A_{ij} b_j;$$

$$a \otimes b = (a_i b_j : 1 \leq i \leq m, 1 \leq j \leq n)$$

With these notational conventions, we define

$$\Theta = \nabla_\varphi - D_t \nabla_{\varphi'} - D_x \nabla_p + D_{tx} \nabla_{p'} + D_{xx} \nabla_q - D_{txx} \nabla_{q'};$$

$$\vartheta = \nabla_{\varphi_\partial} - D_t \nabla_{\varphi_\partial'} + \langle \mathbf{n}(\xi), \nabla_{p_\partial} - D_t \nabla_{p_\partial'} \rangle - \langle \mathbf{n}(\xi) \otimes D_\xi, \nabla_{q_\partial} - D_t \nabla_{q_\partial'} \rangle - D_\xi^* \left[ \nabla_{p_\partial} + \langle \mathbf{n}(\xi) | \nabla_{q_\partial} - D_t \nabla_{q_\partial'} | \right]$$

We have used a hybrid system of notation for these operators, the standard notation of Functional Analysis, and the traditional bracket plus solidus notation of quantum mechanics.

The operators $\Theta_0$, $\Theta_T$, $\vartheta_0$, $\vartheta_T$ are defined by the same formulae but using state variables subscripted with 0, T, $0\partial$, $T\partial$. (Some subscripted variables may be absent from $H$, and then the partial derivatives with respect to the missing terms vanish.) The Hamilton-Euler-Lagrange loaded integro-PDEs are

$$[\psi \; \psi_0 \; \psi_T \; \omega \; \omega_0 \; \omega_T] = [\Theta \; \Theta_0 \; \Theta_T \; \vartheta \; \vartheta_0 \; \vartheta_T] H;$$

$$[\psi \; \omega]|_{t=0^+} = ([\Theta \; \vartheta] H)|_{t=0+}, \quad [\psi \; \omega]|_{t=T^-} = ([\Theta \; \vartheta] H)|_{t=T^-}$$

An extremality principle of the Pontryagin – Boltyanskiy type also holds. We omit both the formulation and the proof of an extremality principle.

## 4. Proofs of selected terms of the Hamiltonian equations

We hall present the proofs for some of the terms of the Hamiltonian equations, from which the proofs of the remaining terms can be inferred. The general approach for discovering the Hamiltonian equations is to introduce penalty terms weighted by the co-states, and thus formulate the Lagrangian functional

$$L = J + \int_0^T \int_\Omega \psi(t,x)[-\varphi(t,x) + RHS(\varphi(t,x))] dx \, dt + \int_0^T \int_\Gamma \omega(t,\xi)[-\varphi(t,\xi) + RHS(\varphi(t,\xi))] d\sigma(\xi) \, dt +$$

$$+ \int_\Omega \left( \psi_0(x)[-\varphi_0(x) + RHS(\varphi_0(x))] + \psi_T(x)[-\varphi_T(x) + RHS(\varphi_T(x))] \right) dx +$$

$$+ \int_\Gamma \left( \omega_0(x)[-\varphi_0(\xi) + RHS(\varphi_0(\xi))] + \omega_T(\xi)[-\varphi_T(\xi) + RHS(\varphi_T(\xi))] \right) d\sigma(\xi)$$

where *RHS* of a component of the state stands for the right-hand side of the corresponding state equation. We then formulate the equations of partial variations of the Lagrangian, i.e. variations with respect to the state only, which, after a conglomeration of calculations, yield the Hamiltonian equations.

We examine the term $J_1 := \int_0^T \int_\Omega F_1(t,x,(\varphi,p,q,\varphi',p',q',u)(t,x)) \, dx \, dt$. We use the symbol $\tilde{\delta}$ for the partial variation with respect to the state only. We have

$$\tilde{\delta} J_1 = \int_0^T \int_\Omega [\nabla_\varphi F_1(t,x,\ldots) \delta\varphi(t,x) + \nabla_p F_1(t,x,\ldots) \delta p(t,x) + \nabla_q F_1(t,x,\ldots) \delta q(t,x) +$$

$$+ \nabla_{\varphi'} F_1(t,x,\ldots) \delta\varphi'(t,x) + \nabla_{p'} F_1(t,x,\ldots) \delta p'(t,x) + \nabla_{q'} F_1(t,x,\ldots) \delta q'(t,x)] dx \, dt$$

We have

$$\int_0^T \int_\Omega \nabla_p F_1(t,x,\ldots) \delta p(t,x) \, dx \, dt = \int_0^T \int_\Gamma \mathbf{n}(\xi) \left( \nabla_p F_1(t,\xi,\ldots) \right) \delta\varphi(t,\xi) \, d\sigma(\xi) \, dt - \int_0^T \int_\Omega \left( D_x \nabla_p F_1(t,x,\ldots) \right) \delta\varphi(t,x) \, dx \, dt;$$

$$\int_0^T \int_\Omega \nabla_{p'} F_1(t,x,\ldots) \delta p'(t,x) \, dx \, dt = \int_\Omega [\nabla_{p'} F_1(T^-,x,\ldots) \delta p(T^-,x) - \nabla_{p'} F_1(0^+,x,\ldots) \delta p(0^+,x)] \, dx +$$

$$+ (-1) \int_0^T \int_\Omega \left( D_t \nabla_{p'} F_1(t,x,\ldots) \right) \delta p(t,x) \, dx \, dt = \int_\Gamma \mathbf{n}(\xi) [\nabla_{p'} F_1(T^-,\xi,\ldots) \delta\varphi(T^-,x) - \nabla_{p'} F_1(0^+,\xi,\ldots) \delta\varphi(0^+,\xi)] \, d\sigma(\xi) +$$

$$+ (-1) \int_\Omega D_x [\nabla_{p'} F_1(T^-,x,\ldots) \delta\varphi(T^-,x) - \nabla_{p'} F_1(0^+,x,\ldots) \delta\varphi(0^+,x)] \, dx - \int_0^T \int_\Gamma \mathbf{n}(\xi) \left( D_t \nabla_{p'} F_1(t,\xi,\ldots) \right) \delta\varphi(t,\xi) \, d\sigma(\xi) \, dt +$$

$$+ \int_0^T \int_\Omega \left( D_{tx} \nabla_{p'} F_1(t,x,\ldots) \right) \delta\varphi(t,x) \, dx \, dt$$

The vanishing of the sum of these variations yields the corresponding terms in the Hamiltonian equations.

## 5 Modelling considerations

It is not feasible to include full modelling examples in the present work. We shall present the axiomatization of Barenblatt models and integral variants of Cattaneo-Vernotte models. The state $\varphi$ will be the transported and diffused quantity (solute concentration in groundwater models, temperature in heat conduction problems). Conservation (of mass or energy) is expressed by

$$\frac{\partial \varphi}{\partial t} + w(t,x,\varphi,\nabla\varphi) = -\nabla\cdot\mathbf{J} + f(t,x,\varphi,\nabla\varphi)$$

where $w$ is a convection term, $\mathbf{J}$ a flux, $f$ a source or sink term. The axiom of Fourier (in heat conduction; the same axiom carries different appellations in the context of fluid mechanics in porous media) is

$$\mathbf{J} = -K(t,x)\nabla\varphi$$

($K$ is a conductivity matrix) thus resulting in a parabolic PDE

$$\frac{\partial \varphi}{\partial t} + w(t,x,\varphi,\nabla\varphi) = \nabla\cdot(K\nabla\varphi) + f(t,x,\varphi,\nabla\varphi)$$

Our axiomatization of a Barenblatt model is achieved by replacing Fourier's axiom with

$$\mathbf{J} = -K(t,x)\nabla\varphi - L(t,x)\nabla\varphi'$$

which leads to

$$\frac{\partial \varphi}{\partial t} + w(t,x,\varphi,\nabla\varphi) = \nabla\cdot(K\nabla\varphi + L\nabla\varphi') + f(t,x,\varphi,\nabla\varphi)$$

The Cattaneo-Vernotte [19, 20] axiom, and its Barenblatt modification (our terminology), are

$$\mathbf{J}' + A(t)\mathbf{J} = -K\nabla\varphi, \quad \mathbf{J}' + A(t)\mathbf{J} = -K\nabla\varphi - L\nabla\varphi'$$

When the retardation factor $A$ is a scalar-valued function of only $t$, the two axioms above lead to hyperbolic PDEs and Barenblatt-augmented hyperbolic PDEs. ("Hyperbolic diffusion" may sound oxymoronic, but it is recognized that the physical concept of diffusion is not identical with second-order PDEs of parabolic type. Furthermore, Hadamard's original classification of second-order linear PDEs does not simply carry over to integro-PDEs.) For a general matricial retardation factor (a matrix-valued function of time and space), the integrated variants of Cattaneo-Vernotte and Cattaneo-Vernotte-Barenblatt are appropriate,

$$\mathbf{J} = \mathbf{J}_0 - \int_0^t K(t,x,s)\nabla\varphi(s,x)\,ds, \quad \mathbf{J} = \mathbf{J}_0 - \int_0^t [K(t,x,s)\nabla\varphi(s,x) + L(t,x,s)\nabla\varphi'(s,x)]\,ds$$

and the last two axioms give rise to (quasi-linear, for this simplified illustration) integro-PDEs, which might be termed <u>integral Cattaneo-Vernotte with Barenblatt terms</u>. If we set $v_0 = -\nabla\cdot\mathbf{J}_0$, the last-mentioned models have the form

$$\frac{\partial \varphi}{\partial t} + w(t,x,\varphi,\nabla \varphi) + v_0 = f(t,x,\varphi,\nabla \varphi) + \nabla \bullet \int_0^t [K(t,x,s)\nabla \varphi(s,x) + L(t,x,s)\nabla \varphi'(s,x)]ds$$

Specifically for the modelling of forest fires, a minimal model (with $\varphi$ interpreted as temperature) includes the integral variant of the Cattaneo-Vernotte axiom, Barenblatt terms, and at least a Fredholm integral term for radiation heat propagation, thus

$$\frac{\partial \varphi}{\partial t} + w(t,x,\varphi,\nabla \varphi) + v_0 = f(t,x,\varphi,\nabla \varphi) + \nabla \bullet \int_0^t [K(t,x,s)\nabla \varphi(s,x) + L(t,x,s)\nabla \varphi'(s,x)]ds +$$
$$+ \int_\Omega R(t,x,y,\varphi(t,y),\nabla_y \varphi(t,y))\,dy$$

Controls are included in *w* and *f*, as well as in the boundary conditions (omitted from this brief discussion).

## 6 Conclusions

We have formulated nonlinear loaded integro-PDEs of relevant types, and we have obtained part of the necessary conditions for optimality (for the concomitant optimal control problems) in the form of what we have termed Hamilton-Euler-Lagrange loaded integro-PDEs.